\documentclass[twocolumn]{autart} 
\usepackage[pdftex,bookmarks,colorlinks]{hyperref}
\usepackage[pdftex,final]{graphicx}
\usepackage{amsmath,amssymb,amsfonts,dsfont}
\usepackage{mathrsfs} 

\begin{document}
\begin{frontmatter}
\runtitle{Minimax estimates for linear DAEs}
\title{Minimax state estimation for linear discrete-time differential-algebraic equations}
\thanks[footnoteinfo]{This paper was presented at 
IFAC Workshop CAO09, Jyvaskyla, Finland, May 6-9, 2009. Corresponding author S. Zhuk. Tel. +38050-52-59138. Fax  +38044-57-56684. }
\author{Sergiy M.Zhuk}
\address{Department of System Analysis and Deicision Making Theory, Taras Shevchenko National University of Kyiv, Ukraine} 
\ead{beetle@unicyb.kiev.ua}
\begin{abstract}
  This paper presents a state estimation approach for an uncertain linear equation with a non-invertible operator in Hilbert space. The approach addresses linear equations with uncertain deterministic input and noise in the measurements, which belong to a given convex closed bounded set. 
A new notion of a\emph{ minimax observable subspace} is introduced. %
By means of the presented approach, new equations describing the dynamics of a minimax recursive estimator for discrete-time non-causal differential-algebraic equations (DAEs) are presented. For the case of regular DAEs it is proved that the estimator's equation coincides with the equation describing the seminal Kalman filter. The properties of the estimator are illustrated by a numerical example.
\end{abstract}
\end{frontmatter}

\section{Introduction}
The importance of models described by DAEs (or descriptor systems) in economics, demography, mechanics and engineering is well known~\cite{Lewis1986}. 
Here, motivated by further applications to linear DAEs, we present a state estimation approach 
for linear deterministic models described by an abstract linear equation in a Hilbert space. 
Our approach is based on ideas underlying $H_2/H_\infty$
filtering~\cite{Balakrishnan1976,Basar1995} and set-membership state estimation~\cite{Bertsekas1971,Tempo1985,Nakonechnii1978,Kurzhanski1997,Milanese1991}.\\ 
$H_2$-estimators like Kalman or
Wiener filters~\cite{Balakrishnan1976,Albert1972} give estimations of the system state with minimum error variance. The $H_2$-estimation problem for linear time-variant DAEs was studied in
\cite{nikh1} without restricting the DAE's matrices. The resulting algorithm requires the calculation of the so-called ``3-block matrix pseudoinverse''. In \cite{ishixara2} the authors introduced explicit formulas for the 3-block matrix pseudoinverse and derived a recurrence filter, assuming a special structure for the DAE's matrices. %
A brief overview of steady-state $H_2$-estimators is presented
in~\cite{Deng1999}.
\newline $H_\infty$ estimators minimize a norm 
of the operator mapping unknown disturbances with finite energy to filtered
errors~\cite{Basar1995}. We stress that the $H_\infty$ estimator coincides with a certain Krein space $H_2$ filter \cite{Sayed2001}. 
The $H_\infty$ filtering technique was applied to linear time-invariant DAEs with regular matrix pencils in~\cite{Xu2007}. \\
A basic notion in the theory of set-membership state estimation is that of an a posteriori set or informational set. This notion has roots in control theory~\cite{Balakrishnan1976}. 
By definition, it is the set of all
possible state vectors $\varphi$, that are consistent with a measured output $y$, 
provided that an uncertain input $f$ and measurement error $\eta$
belong to some bounded set $\mathscr G$. %
We will be interested in the case when the state $\varphi\in\mathcal H$ obeys an abstract linear
equation $L\varphi=f$, provided $y=H\varphi+\eta$, $(f,\eta)\in \mathscr
G$, where $\mathscr G$ is a bounded closed convex subset of an abstract
Hilbert space. %
The problem is to find an estimation $\hat\varphi$ of $\varphi$ with minimal
worst-case error. This problem was previously considered in~\cite{Nakonechnii1978,Tempo1985}. Due to~\cite{Nakonechnii1978} a vector $
\hat\varphi$ is called a \emph{linear minimax
  a-posteriori estimation} (or a \emph{central
  algorithm} due to~\cite{Tempo1985}) iff $\forall \ell\in\mathcal H$
\begin{equation}
  \label{eq:mcenter}
  \langle\ell,\hat\varphi\rangle=(
 \sup_{\mathscr G(y)}\langle\ell,\varphi\rangle+
 \inf_{\mathscr G(y)}\langle\ell,\varphi\rangle)/2, 
\end{equation}
provided that an a posteriori set $\mathscr
G(y):=\{\varphi:(L\varphi,y-H\varphi)\in\mathscr G\}$ is a bounded
convex subset of the Hilbert space $\mathcal H$. 
Note that if there exists $\varphi_0$ so that $L\varphi_0=0$, $H\varphi_0=0$, then $\sup_{\mathscr  G(y)}\langle\ell,\varphi\rangle=+\infty$ for some $
\ell$. Thus, the above approach does not work if $L$ is non-injective. In this paper we generalize the approach of~\cite{Tempo1985,Nakonechnii1978} to linear equations with non-injective $L$. Futher generalization is presented in~\cite{Zhuk2009d}.  

The main contribution of this paper is a new notion of a \emph{minimax
observable subspace} $ 
\mathcal L
$ for the pair $(L,H)$ (Definition~\ref{d:1}). 
It is useful when one needs to evaluate \emph{a priori} how far the estimation $\hat\varphi$ is from a ``real'' state $\varphi$ in the worst case, provided $\hat \varphi$ is constructed from the measurements $y$. Due to Proposition~\ref{p:gnrl}, the \emph{worst-case estimation error} is finite iff $\mathcal L=\mathcal H$; otherwise $\hat\varphi$ may be too far from a ``real'' state $\varphi$ for some directions $\ell\in\mathcal H$, 
even for bounded $f$ and $\eta$. In fact, given $y$, we can provide an 
estimation with finite worst-case error for the projection of $\varphi$ onto $\mathcal L$ only. Thus $\mathcal L$ is an analog of the observable subspace~\cite[p.240]{Balakrishnan1976} for the pair $(L,H)$ 
in the context of set-membership state estimation. \\
The introduced notion allows the generalization of ideas from~\cite{Tempo1985,Nakonechnii1978,Tempo1988} to non-injective linear mappings, in particular for the case $
 L\varphi(t)=((F \varphi)_t-C(t)\varphi(t),F\varphi(t_0))$ with $F\in\mathbb R^{m\times n}$ which arise in the state estimation for linear continuous non-causal\footnote{DAE is said to be non-causal if the corresponding
    initial-value problem has more than one solution. } DAEs~\cite{Zhuk2007}.
As a consequence, one can apply the minimax framework, originally developed \cite{Bertsekas1971,Baras1995} for DAEs ($F=E$ in the linear case) with bounded uncertainties, to DAEs~\cite{Zhuk2009d} with unbounded inputs (see example in Section 3).\\
In order to stress connections with $H_\infty$ approach, we note that the minimax framework \cite{Bertsekas1971,Baras1995} incorporates the set-membership state estimation and $H_\infty$ filtering for ODEs by application of dynamic programming~\cite{Bellman1972} to the informational state $X(\tau)$: for linear ODE the worst-case estimation is
set to be the Tchebysheff center of $X(\tau)$.  
Although we derive the estimation from the minimization of the worst-case error, as it is stated in Definition 1, our approach (for the ellipsoidal bounding set
$\mathscr G$ and causal\footnote{Note that the dynamic programming was previously applied to causal DAEs in~\cite{Bender1987} in order to construct a regulator in LQ-control problems with DAE constraints.} DAEs) results in the same estimation and error as in~\cite{Bertsekas1971}. Thus, the $\ell$-minimax estimation gives a proper generalization of the recurrence algorithm from \cite{Bertsekas1971} to the case when $X(\tau)$ may be unbounded.\\
We illustrate the benefits of the new notion by introducing a minimax recursive estimator for discrete-time non-causal DAEs: 
it works for non-causal DAEs unlike~\cite{ishixara2,Deng1999,Xu2007,Zhang2003} and for the regular case it coincides with one proposed in~\cite{ishixara2} (Corollary~\ref{l:eqfltr}). In addition, the minimax observability subspace allows one to identify the observable (in the minimax sense) part of the state with respect to given measurements. Computing the index of non-causality, one 
can a priori check how good connections between observations and state are: models with zero index are fully observable while models with non-zero index have an unobservable part in the state. \\
This paper is organized as follows. In Section 2 we give definitions (Definition 1) of the minimax estimation, error and observable subspace for abstract linear equations and we construct the estimation for a convex bounded $\mathscr G$, in particular for an ellipsoidal $\mathscr G$ (Proposition 1). In Section 3 we introduce the minimax observable subspace and index of non-causality for DAEs in discrete time (Definition 2) and we derive the minimax estimator (Theorem 1). Also we discuss connections to $H_2/H_\infty$ framework (Corollary 1) and present an example. 
\newline{\bf Notation}. \emph{Linear mappings:} 
$\langle\cdot,\cdot\rangle$ denotes the inner product; $
\mathscr L(\mathcal H_1,\mathcal H_2)$ denotes  the space 
of all bounded linear mappings from $
\mathcal H_1$ to $\mathcal H_2$, 
$\mathscr L(\mathcal H):=\mathscr L(\mathcal H,\mathcal H)$; 
$\mathds{1}_{\mathcal H}$ is the identity mapping in $\mathcal H$; 
$\mathscr{D}(L)$, $R(L)$, $N(L)$ denote, respectively, the domain, range, and null-space of a linear mapping $L:\mathscr D(L)\mapsto R(L)$; 
$L^*:\mathcal H\to\mathcal H$ is the adjoint of $L$;  
$F'$ denotes the transpose of $F$; $F^+$ is 
the pseudoinverse of $F$; $E$ is the identity matrix;  
$\mathrm{diag}(A_1\dots A_n)$ denotes a diagonal matrix with $A_i$, $i=\overline{1,n}$ on its diagonal; $\{x_s\}_1^n:=(x_1,\dots,x_n)$ is an element of $\mathcal H_1\times\dots \times \mathcal H_n$, $0_{mn}\in\mathbb R^{m\times n}$ denotes the $m\times n$-zero matrix.
\newline\emph{Functionals:}
$
I_1(x):=\langle Q_1 Lx,Lx \rangle+\langle Q_2 Hx,Hx\rangle
$, $
I(x):=\langle Q_1Lx,Lx\rangle+\langle Q_2(y-Hx),y-Hx\rangle$ with positive definite self-adjoint 
$Q_1\in\mathscr L(\mathcal F)$ and $Q_2\in\mathscr L(\mathcal Y)$;
$c(G,x):=\sup\{\langle x,y\rangle,y\in G\}$ is a support function of $G$; 
$\gamma_\pm:=\frac 12(c(\mathscr G(y),\ell)\pm c(\mathscr
G(y),-\ell))$; 
$\|x\|^2_{S}=\langle S x,x\rangle$.
\newline\emph{Sets:}
$
\mathscr G^\beta(0):=\{x: I_1(x)\le\beta\}$, 
$\mathscr G(0):=\{x: I_1(x)\le1\}$; $
\mathrm{dom} f:=\{x:f(x)<\infty\}$ is an effective domain of $f$;
$\mathrm{Argmin}_xf:=\{x:f(x)=\min_x f\}$ is the set of global minima of $f$;
$\overline{G}$ is the closure of a set $G$.
\section{Linear minimax estimation problem in a Hilbert space}
\label{s:mnmx}
Let vector $y$ be observed in the form of 
\begin{equation}
  \label{eq:y}
  y=H\varphi+\eta
\end{equation}
where $\varphi$ obeys the equation 
\begin{equation}
  \label{eq:Lfi}
  L\varphi=f
\end{equation}
We assume that $
L:\mathscr D(L)\mapsto \mathcal F$ is a closed linear operator~\cite[p.63]{Balakrishnan1976}, $\mathscr D(L)\subset\mathcal H$ is a linear set, $\overline{\mathscr{D}(L)}=\mathcal{H}$, $H\in\mathscr L(\mathcal H,\mathcal Y)$. Also we assume that $(f,\eta)$ is an unknown element of a convex bounded closed set $\mathscr G\subset\mathcal F\times\mathcal Y$, $\mathcal H,\mathcal F,\mathcal Y$ are Hilbert spaces.
\begin{defn}\label{d:1}
Let  $\mathcal L:=\{\ell\in\mathcal H:\hat\rho(\ell)<\infty\}$ with $$
\hat\rho(\ell):=\inf_{\varphi\in\mathscr G(y)}\rho(\ell,\varphi),\quad
\rho(\ell,\varphi):=\sup_{\psi\in\mathscr G(y)}|\langle \ell,\varphi-\psi\rangle|
$$ The set $\mathcal L$
is called a minimax observable subspace for the pair $(L,H)$. A vector
$\hat\varphi\in\mathscr G(y)$ is called a minimax a posteriori estimation in
the direction\footnote{The subspace $l:=\{\alpha\ell:\alpha\in\mathbb R\}$, assigned with $\ell\in\mathcal{H}$, defines some 
  direction in $\mathcal H$. 
To estimate $\varphi$ in the direction $\ell$ means to estimate the projection $\langle\ell,\varphi\rangle\ell$ of $\varphi$ onto $l$. 
} $\ell$ ($\ell$-estimation) if $\rho(\ell,\hat\varphi)=\hat\rho(\ell)$. The number $\hat\rho(\ell)$ is called a
minimax a posteriori error in the direction $\ell$ ($\ell$-error).  
\end{defn}
\emph{Our aim here is}, given $y$, to construct the $\ell$-estimation
  $\hat\varphi$ of the state $\varphi$, $\ell$-error
$\hat\rho(\ell)$ and minimax observable subspace $\mathcal L$, provided $
\ell\in\mathcal L$. Note that $\hat\rho(\ell)=+\infty$ if $\ell\not\in\mathcal L$ so that any $\psi\in\mathcal H$ is a $\ell$-estimation by Definition~\ref{d:1}.  
\begin{prop}\label{p:gnrl}
  Let $\mathscr G$ be a convex closed bounded subset of $
\mathcal F\times\mathcal Y$. Then $\ell\in\mathcal L$ iff  
$  \ell,-\ell\in
  \mathrm{dom}\,c(\mathscr G(y),\cdot)
$. If $\ell\in\mathcal L$ then the $\ell$-estimation $\hat\varphi$ 
along with $\ell$-error $\hat\rho(\ell)$ obey
\begin{equation}
  \label{eq:hlfi}
\langle \ell,\hat\varphi\rangle=\gamma_-,\hat\rho(\ell)=\gamma_+
\end{equation} 
Define $T:\mathscr D(T)\to\mathcal H$ by the rule $T(z,u):=L^*z+H^*u$ with $
\mathscr D(T):=\mathscr D(L^*)\times\mathcal Y$ and let $$
\mathscr G=\{(f,\eta):\langle Q_1f,f\rangle+\langle Q_2\eta,\eta\rangle\le 1\}
$$ If $R(T)=\overline {R(T)}$ then $\hat
x\in\mathrm{Argmin}_x\,I$ is the $\ell$-estimation, 
$\mathcal L=\mathrm{dom}\,c(\mathscr G(0),\cdot)= R(T)$ and 
  \begin{equation}
    \label{eq:qeser}
  \hat\rho(\ell)=
(1-I(\hat x))^\frac 12 c(\mathscr G(0),\ell).
 \end{equation} 
The worst-case estimation error for any direction is 
\begin{equation}
 \label{eq:c:vecest}
 \begin{split}
    & \sup_{x\in\mathscr G(y)}\|\hat x-x\|=
    \inf_{\varphi\in\mathscr G(y)}\sup_{x\in\mathscr G(y)}\|\varphi-x\|=\\
    &(1-I(\hat x))^\frac 12 
    \sup_{\|\ell\|=1}
    c(\mathscr G(0),\ell)=
    \sup_{\|\ell\|=1}\hat\rho(\ell)<+\infty
  \end{split}
\end{equation}
If $\mathcal L=\mathcal H$ then~\eqref{eq:c:vecest} is finite. 
\end{prop}
\begin{pf}
Let $\ell(\mathscr G(y))=\{\langle\ell,\psi\rangle,\psi\in\mathscr G(y)\}$. 
Since $\mathscr G(y)$ is convex (due to convexity of $\mathscr D(L)$ and $\mathscr G$) and $x\mapsto\langle\ell,x\rangle$ is continuous, it follows that $\overline{\ell(\mathscr G(y))}$ is connected. Noting that 
$\inf_{\mathscr G(y)}\langle\ell,\psi\rangle=-c(\mathscr G(y),-\ell)$, 
we see $$
\overline{\ell(\mathscr G(y))}=
[-c(\mathscr G(y),-\ell),c(\mathscr G(y),\ell)]
\subset\mathbb R^1
$$ Thus $\rho(\ell,\varphi)=+\infty$ if $\ell,-\ell\not\in \mathrm{dom}\,c(\mathscr
G(y),\cdot)$. Otherwise $\overline{\ell(\mathscr G(y))}$ is bounded, implying $
\ell\in\mathcal L$. Hence,  $\ell\in\mathcal L$ iff $\ell,-\ell\in
  \mathrm{dom}\,c(\mathscr G(y),\cdot)$.\\
Let  $
\ell\in\mathcal L$, $\varphi\in\mathscr G(y)$. Since $\overline{\ell(\mathscr G(y))}$ is connected, there exists 
$\varphi_*\in\mathscr G(y)$ so that $\langle\ell,\varphi_*\rangle=\gamma_-$ is the central point of $\overline{\ell(\mathscr G(y))}$. 
The worst-case distance $\rho(\ell,\varphi)$ is equal to the sum of the distance $|\langle\ell,\varphi-\varphi^*\rangle|$ between $\langle\ell,\varphi\rangle$ and the central point  $\gamma_-$ and the distance $\gamma_+$ between one of the boundary points of $\overline{\ell(\mathscr G(y))}$ and  $\gamma_-$. Therefore, $\gamma_-$ has the minimal worst-case distance $\gamma_+$. Hence, $\varphi_*=\hat\varphi$ due to Definition~\ref{d:1}, which implies~\eqref{eq:hlfi}. \\
We proceed with the ellipsoidal $\mathscr G$. Let $Q_1=\mathds{1}_{\mathcal F},Q_2=\mathds{1}_{\mathcal Y}$ for a simplicity. 
Due to~\cite[Sec 5.\S 3]{Kato1966} $R(T)=\overline{R(T)}$ implies $
R(T^*)=\overline{R(T^*)}$. Thus~\cite[p.14,Cor.1.4.3]{Balakrishnan1976}, there exists $\hat x\in\mathscr D(L)$ so that $T^*\hat x$ is the projection of $
(0,y)$ onto $R(T^*)=\{(Lx,Hx),x\in\mathscr D(L)\}$, implying $\hat x\in\mathrm{Argmin }_xI$, and $$
\langle y-H\hat x,Hx\rangle=\langle L\hat x,Lx\rangle, \forall x\in\mathscr D(L)$$Noting this, one easily derives\footnote{$I(\hat x-x)=I_1(x)+I(\hat x)-2\langle L\hat x,L x\rangle+2\langle y-H\hat x,Hx\rangle$} $
I(\hat x- x)=I_1(x)+ I(\hat x)
$ for all $x\in\mathscr D(L)$. Having it in mind and noting that $
\mathscr G(y)=\{\varphi:I(\varphi)\leqslant1\}$ and $\hat x\in\mathscr G(y)$, 
one derives 
\begin{equation}
  \label{eq:GyG0}
\hat x+\mathscr G^\beta(0)= \mathscr G(y)  
\end{equation}
where $\beta:=1-I(\hat x)$. \eqref{eq:GyG0} 
implies~\cite[p.113]{Rockafellar1970}  
\begin{equation}
  \label{eq:cGy}
  c(\mathscr G(y),\ell)=\langle \ell,\hat x\rangle+c(\mathscr G^\beta(0),\ell),
\forall\ell\in\mathcal H
\end{equation} 
The definition of $\gamma_-$, \eqref{eq:cGy} and $c(\mathscr G^\beta(0),\ell)=c(\mathscr
G^\beta(0),-\ell)$ imply $
\gamma_-=\langle \ell,\hat x\rangle$. Due to \eqref{eq:hlfi}, $\hat x$ is the
$\ell$-estimation. 
\newline Let us prove~\eqref{eq:qeser}. 
If $x\in\mathscr G^\beta(0)$ then $I_1(\beta^{-\frac12}x)\le 1$. Thus $\beta^{-\frac12}\mathscr G^\beta(0)\subset\mathscr G(0)$ implying $
\mathscr G^\beta(0)\subset\beta^\frac 12\mathscr G(0)$. 
If $
x\in\mathscr G(0)$ then $I_1(\beta ^{\frac 12}x)=\beta I_1(x)\le\beta\Rightarrow\beta^\frac 12\mathscr G(0)\subset\mathscr G^\beta(0)$. 
Therefore
\begin{equation}
  \label{eq:GbG}
  \mathscr G^\beta(0)=\beta^{\frac 12}\mathscr G(0)\Rightarrow
c(\mathscr G^\beta(0),\ell)=\beta^{\frac 12}c(\mathscr G(0),\ell)
\end{equation} 
Now \eqref{eq:hlfi} and \eqref{eq:cGy} imply~\eqref{eq:qeser}. Hence,  $
\mathcal L=\mathrm{dom}\,c(\mathscr
G(0),\cdot)$ due to Definition~\ref{d:1}. 
\newline Let us prove $R(T)=\mathcal L$. 
If $\ell\in R(T)$ then $
\ell=L^*z+H^*u$ for some $z\in\mathscr D(L^*)$, 
$u\in\mathcal Y$ and we get $
\forall\varphi\in\mathscr G(0)$ by Cauchy inequality~\cite[p.4]{Balakrishnan1976} $$
\langle\ell,\varphi\rangle=\langle z,L\varphi\rangle+
\langle u,H\varphi\rangle\le\|z\|^2+\|u\|^2<+\infty
$$  so that $R(T)\subset\mathrm{dom}\,c(\mathscr
G(0),\cdot)=\mathcal L$. 
If $\ell\notin R(T)$ then $\langle\ell,x\rangle>0$ for some $x\in N(T^*)$ as $\mathcal H = \overline{R(T)}\oplus N(T^*)$. Noting that $\mathscr G(0)=\{x:\|T^*x\|^2\le1\}$ we derive $
c(\mathscr G(0),\ell)\ge\sup\{\langle\ell,x\rangle:T^*x=0\}=+\infty
$. Hence, $\mathcal L\subset R(T)$. 
\newline Let us prove~\eqref{eq:c:vecest}. Set $\hat\alpha:=\inf_{\varphi\in\mathscr G(y)}\sup_{x\in\mathscr G(y)}\|\varphi-x\|$. 
Using Definition~\ref{d:1}, one derives\footnote{Note~\cite[p.42]{Balakrishnan1976} that  
$\|\varphi\|=\sup_{\|\ell\|=1}\langle\ell,\varphi\rangle$ and~\cite[p.55]{Balakrishnan1976} $
\inf_x\sup_yF(x,y)\ge\sup_y\inf_xF(x,y)$ for 
convex-concave $F$.} 
\begin{equation*}
  \begin{split}
    &\hat\alpha=
    \inf_{\varphi\in\mathscr G(y)}\sup_{x\in\mathscr G(y)}
    \sup_{\|\ell\|=1}|\langle\ell,\varphi-x\rangle|\ge\\
    &\sup_{\|\ell\|=1}\inf_{\varphi\in\mathscr G(y)}
    \sup_{x\in\mathscr G(y)}|\langle\ell,\varphi-x\rangle|=
    \sup_{\|\ell\|=1}\hat\rho(\ell)=\\
    &\sup_{\|\ell\|=1}\sup_{x\in\mathscr
      G(y)}|\langle\ell,\hat x-x\rangle|=
    \sup_{x\in\mathscr G(y)}\|\hat x-x\|\ge\hat\alpha
  \end{split}
\end{equation*}
Now assume $\mathcal L=\mathcal H$. 
Since $\mathrm{dom}\,c(\mathscr G(0),\cdot)=
\mathcal L$, it follows $c(\mathscr G(0),\cdot)$ is finite in $\mathcal H$ and therefore~\cite[\S 2.3]{Ekeland1976} continuous. As a consequence, \eqref{eq:c:vecest} is finite. 
\end{pf}
\section{ $\ell$-estimation for non-causal DAEs }
\label{s:obs}
\begin{figure}[h]\centering
\begin{minipage}[c]{550pt} 
\includegraphics
[width=8.4cm,height=5cm]
{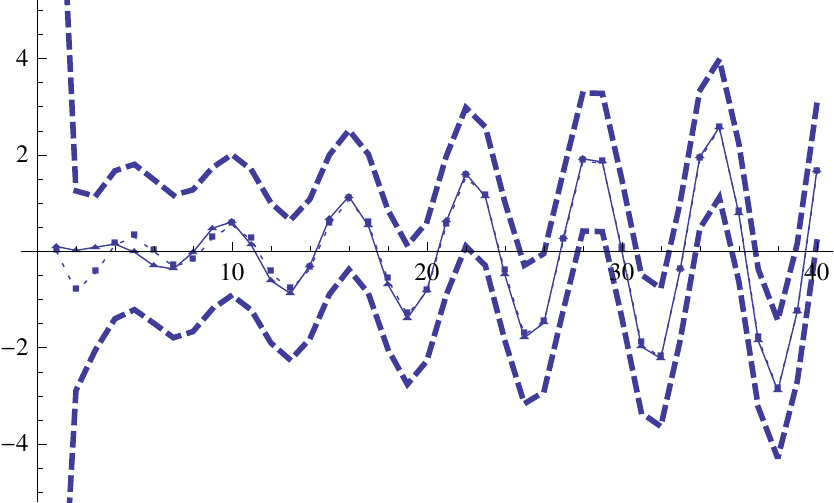}
\end{minipage}
\caption{The dashed line corresponds to the real values of $\langle\ell,p^*_k\rangle$ with $\ell=(0,1)$, $k=\overline{1,40}$; 
the solid line corresponds to the $\ell$-estimation $\langle\ell,\hat x_k\rangle$; 
the bold dashed lines represent dynamics of the boundary points of the segment $\ell(\tilde X(\tau))$. 
Note that the trajectory of the estimation is centered with respect to ``the bounds'' -- bold dashed lines.}
\label{fig:2}
\end{figure}

Consider the model
\begin{align}
    &F_{k+1}x_{k+1}-C_k x_k = f_{k+1}, F_0x_0=f_0,\label{eq:state}\\
    &y_k=H_k x_k+g_k, k=0,1,\dots\label{eq:mes}
\end{align}
where $F_k,C_k\in\mathbb R^{m\times n}$, $H_k\in
\mathbb R^{p\times n}$, $x_k\in\mathbb R^n$ is a state, $f_k\in\mathbb R^m$ is an input and $y_k,g_k\in \mathbb R^p$ represent an output and
the output's noise respectively. 
In what follows we assume that an initial state $x_0$ belongs to
the affine set $\{x:F_0x=f_0\}$. We define 
$\xi_\tau=(\{f_s\}_0^\tau,\{g_s\}_0^\tau)$ 
and assume 
\begin{equation}
  \label{eq:ellips}
\xi_\tau
\in\mathscr G=
\{\xi_{\tau}:
\sum_0^{\tau}\langle S_i f_i,f_i\rangle+\langle R_ig_i,g_i\rangle
\leqslant1 \}
\end{equation}
where $S_k\in\mathbb R^{m\times m}$ and $R_k\in \mathbb R^{p\times p}$ are positive definite self-adjoint matrices.\\
Suppose we observe $y^*_1,\dots, y^*_\tau$, provided that $y_k^*$ is derived from~\eqref{eq:mes} with $
 g_k=g^*_k$ and $x_k=x^*_k$, which obeys~\eqref{eq:state} with $f_k=f^*_k$, and 
 $(\{f^*_s\}_0^\tau,\{g^*_s\}_0^\tau)\in\mathscr G$. 
Denote by $X(\tau)$ the set of all possible states $x_\tau$ of~\eqref{eq:state} consistent with measurements $y^*_1,\dots, y^*_\tau$ and uncertainty description~\eqref{eq:ellips}. 
\begin{defn}\label{d:2}
We say that $\overline x_\tau$ is a $\ell$-estimation of the state $x^*_\tau$ in the direction $\ell\in\mathbb R^{n}$ iff $$
\sup_{x\in X(\tau)}|\langle\ell,x-\overline x_\tau\rangle|=
\hat{\rho}(\ell,\tau):=\inf_{z\in X(\tau)}\sup_{x\in X(\tau)}|\langle\ell,x-z\rangle|
$$
$\hat{\rho}(\ell,\tau)$ is said to be 
an $\ell$-error. A minimax observable subspace at the
instant $\tau$ for the model~\eqref{eq:state}-\eqref{eq:mes} is defined by $\mathcal L(\tau)=\{\ell:\hat{\rho}(\ell,\tau)<\infty\}$.
$I_\tau:=n-\mathrm{dim}\, \mathcal{L}(\tau)$ is called an index of non-causality of the model~\eqref{eq:state}-\eqref{eq:mes}.
\end{defn}
\begin{thm}\label{t:mnmx}
  Define $\hat\beta_\tau:=1-\alpha_\tau+
\langle P_\tau\hat x_\tau,\hat x_\tau\rangle{}$ and $\hat x_\tau=P_\tau^+r_\tau$ with $r_0=H'_0R_0y_0$, $\alpha_0=\langle R_0y_0,y_0\rangle$, 
\begin{equation*}
    \begin{split}
      &P_k=
      H'_kR_kH_k+F'_k[S_{k-1}-S_{k-1}C_{k-1}B_{k-1}^+C'_{k-1}S_{k-1}]
      F_k,\\
      &P_0=F'_0S_0F_0+H'_0R_0H_0,B_k=P_{k}+C'_{k}S_{k}C_{k},\\
      &\alpha_k=\alpha_{k-1}+\langle R_ky_k,y_k\rangle-\langle B_{k-1}^+
 r_{k-1},r_{k-1} \rangle, \\
 &r_k=F'_kS_{k-1}C_{k-1}B_{k-1}^+r_{k-1}+H'_kR_k y_k,
          \end{split}
  \end{equation*}
Then $\hat x_\tau$ is the 
$\ell$-estimation of $x^*_\tau$, $\mathcal L(\tau)=\{\ell:P_\tau^+P_\tau\ell=\ell\}$, $
    \hat{\rho}(\ell,\tau)={\hat\beta}_\tau^\frac12
    \langle P^+_\tau\ell,\ell\rangle^\frac 12$ and~\footnote{Note that $\tau\mapsto X(\tau)$ represents the a posteriori set-valued observer~\cite{Shamma1999}}
  \begin{equation}
    \label{eq:Xtau}
    X(\tau)=\hat x_{\tau}+\hat\beta_\tau^{\frac 12}\tilde{X}(\tau),
    \tilde X(\tau):=\{x:\langle P_\tau x,x\rangle\le 1\}
 \end{equation}
\end{thm}
\begin{pf}
In order to apply Proposition~\ref{p:gnrl}, we rewrite~\eqref{eq:state}-\eqref{eq:mes} 
in the operator form: set $\mathcal{H}=(\mathbb R^n)^{\tau+1}$, $\mathcal Y=(\mathbb R^p)^{\tau+1}$, $\mathcal F=(\mathbb R^m)^{\tau+1}$ and $\varphi^*=\{x^*_s\}_0^\tau$, $ y^*=\{y^*_s\}_0^\tau$, $
f^*=\{f^*_s\}_0^{\tau}$, $\eta^*=\{g^*_s\}_0^\tau$, $$
L=\left(\begin{smallmatrix}
      F_0&&0_{mn}&&0_{mn}&&\dots&&0_{mn}&&0_{mn}\\
      -C_0&&F_1&&0_{mn}&&\dots&&0_{mn}&&0_{mn}\\
      \dots&&\dots&&\dots&&\dots&&\dots&&\dots\\
      0_{mn}&&0_{mn}&&0_{mn}&&\dots&&-C_{\tau-1}&&F_\tau
\end{smallmatrix}\right),
 H=\left(\begin{smallmatrix}H_0&&\dots&& 0_{pn}\\
 \dots&&\dots&&\dots\\
 0_{pn}&&\dots&& H_\tau\end{smallmatrix}\right)
$$ 
Define $\mathcal
P_\tau = \bigl(\begin{smallmatrix}
  0_{nn},&&\dots,&&0_{nn},E\end{smallmatrix}\bigr)$ and rewrite \eqref{eq:ellips} with $
Q_1=\mathrm{diag}\left(S_0,\dots,S_{\tau}\right)$,
$Q_2=\mathrm{diag}\left(R_0,\dots,R_{\tau}\right)$ as $$ 
\xi_\tau=(f,\eta)\in \mathscr G=
\{(f,\eta):\langle Q_1f,f\rangle+\langle Q_2\eta,\eta\rangle\le 1\}.
$$ It is clear that $
y^*$, $H$, $\eta^*$, $L$, $f^*$, $\varphi^*$, defined as above, satisfy~\eqref{eq:y}-\eqref{eq:Lfi} and $(f^*,\eta^*)\in \mathscr G$. Let $\mathscr G(y^*)$ denote the a posteriori set generated by $y^*$. Then $X(\tau)=\mathcal P_\tau(\mathscr
G (y^*))$ by definition. Thus
\begin{equation}
  \label{eq:rhot}
  \hat{\rho}(\ell,\tau)=\inf_{\varphi\in\mathscr G(y^*)}\sup_{\psi\in\mathscr G(y^*)}|\langle{\ell,\mathcal P_\tau(\varphi-\psi)}|=\hat\rho(l)
\end{equation} with 
$l=\mathcal P_\tau'\ell$. Hence, $\overline x_\tau=\mathcal P_\tau\hat\varphi$, where $\hat\varphi$ is the $l$-minimax estimation of the state $\varphi^*$ of 
\eqref{eq:Lfi} in the sense of Definition~\ref{d:1}. Proposition~\ref{p:gnrl} implies $
\hat\varphi=\hat x$ and $
\hat\rho(l)=\beta^\frac 12 c(\mathscr G(0),l)
$. Let us prove~\eqref{eq:Xtau}. 
\eqref{eq:GyG0} and \eqref{eq:GbG} imply $\mathcal P_\tau(\mathscr
G (y^*))=\mathcal P_\tau\hat x+\beta^\frac12\mathcal P_\tau(\mathscr G(0))$. 
Therefore, $X(\tau)=\mathcal P_\tau\hat x+\beta^\frac 12\mathcal P_\tau(\mathscr G(0))$. 
Now, let us prove $\hat x_\tau=\mathcal P_\tau\hat x$ by the direct calculation.
 Define $$ 
V_\tau(x_0,\dots, x_\tau):=
\Phi(x_0)+\sum_{s=0}^{\tau-1}\Phi_s(x_s,x_{s+1})
$$ with $\Phi_k(x,p):=\|F_{k+1}x-C_{k}p\|_{S_{k+1}}^2+\|y_{k+1}-H_{k+1}x\|_{R_{k+1}}^2$,
$\Phi(x):=\|F_0x\|_{S_0}^2+\|y_0-H_0x\|_{R_0}^2$, and $$
\mathcal B_\tau(p):=\min_{x_0\dots x_{\tau-1}}V_\tau(x_0\dots x_{\tau-1}, p),
\mathcal B_0(p):=\Phi(p)
$$
\textbf{Lemma 1.} 
\emph{Let $p\in\mathbb R^n$. There exists $(\tilde x_1\dots\tilde x_{\tau-1})\in\mathbb( R^n)^{\tau-1}$ so that $V_\tau(\tilde x_1\dots\tilde x_{\tau-1},p)=\mathcal B_\tau (p)$ and 
\begin{equation}
  \label{eq:BPk}
  \mathcal B_k(p)=\langle P_kp,p\rangle-2\langle r_k,p\rangle+\alpha_k,P_k\ge 0  
\end{equation}}
Lemma 1 implies $\mathcal B_\tau{}$ is a quadratic and non-negative function. Therefore $
\mathrm{Argmin }\mathcal B_\tau=\{x:P_\tau x=r_\tau\}\ne\varnothing
$. This and $I(\varphi)=V_\tau(x_0,\dots, x_\tau)$ imply $$
\mathcal B_{\tau}(\hat x_\tau)=\min \mathcal B_\tau=\min V_{\tau}=\min I
$$ Defining $\hat x:=(\tilde x_1\dots\tilde x_{\tau{}-1},\hat x_\tau{})$ with $\tilde x_k$ taken as in Lemma 1 for $p=\hat x_\tau$, we obtain $
I(\hat x)=\min I$ and $\hat x_\tau{}=\mathcal P_\tau{}\hat x$. Therefore, $\hat x_\tau{}$ is a $\ell$-estimation.
\newline We note $\min I=\mathcal B_\tau(\hat x_\tau)=
1-\hat\beta_{\tau}$. Thus $\beta=\hat\beta_{\tau}$ by definition. \\
Let us prove $\mathcal P_\tau(\mathscr G(0))=\tilde X(\tau)$. Since $\mathscr G(0)$ does not depend on $y^*$ and $\mathscr G(0)=\mathscr G(y^*)$ provided $y^*=0$, we can calculate $\mathcal P_\tau(\mathscr G(0))$ assuming $y^*=0$. 
In this case $I=I_1=V_\tau$, $\hat x_\tau=0$ and $\mathcal B_\tau(x)=\langle P_\tau x,x\rangle$ so that $x\in\tilde X(\tau)\Leftrightarrow\mathcal B_\tau(x)\le1$.  
If $x\in \mathcal P_\tau(\mathscr G(0))$ then, by definition, there exist $x_{1}\dots x_{\tau-1}$ so that $V_{\tau}(x_{1}\dots  x_{\tau-1}, x)\le 1$, implying $\mathcal B_{\tau}(x)\le 1$.  
Now, let $\mathcal B_{\tau}(x)\le 1$. Then $V_\tau(\tilde x_1\dots\tilde x_{\tau-1},x)=\mathcal B_\tau (x)$ due to Lemma 1 and thus
$x\in\mathcal P_{\tau}(\mathscr G(0))$ by definition. 
\newline Formulae \eqref{eq:qeser}, \eqref{eq:rhot} and $\mathcal P_\tau(\mathscr G(0))=\tilde X(\tau)$ imply 
$\hat\rho(\ell,\tau)=
\hat\beta^{\frac 12}_\tau c(\tilde X(\tau),\ell)={\hat\beta}_\tau^\frac12
\langle P^+_\tau\ell,\ell\rangle^\frac 12$ for $\ell\in\mathcal{L}(\tau)$ 
and $\mathcal L(\tau)=\mathrm{dom}\,c(\tilde X(\tau),\cdot)=\{\ell:P_\tau^+P_\tau\ell=\ell\}
$. Details of calculation of $c(\tilde X(\tau),\cdot)$ are given in~\cite[p.108]{Rockafellar1970}. This completes the proof. 
\par\textbf{Proof of Lemma 1.}  
We shall apply the dynamic programming~\cite{Bellman1972}. Since $V_\tau$ is additive, it follows that
\begin{equation}
  \label{eq:Bdfn}
\mathcal B_\tau(p)=\min_{x_{\tau-1}}\{\Phi_{\tau-1}(x_{\tau-1},p)+
\mathcal B_{\tau-1}(x_{\tau-1})\}
\end{equation}
$V_\tau$ is convex and non-negative by definition. 
Thus $\mathcal B_\tau$ is non-negative and convex for any $\tau\in\mathbb N$. 
Convexity is implied by the definition of $\mathcal B_\tau$ as 
 for any convex function $ 
(x,y)\mapsto f(x,y)$ the function $
y\mapsto\min_xf(x,y)$ is convex~\cite[p.38]{Rockafellar1970}.

Let us prove \eqref{eq:BPk} by induction. 
\eqref{eq:BPk} holds for $\mathcal B_0$ and $P_0$. We shall derive \eqref{eq:BPk} for $\mathcal B_{k+1}$, $P_{k+1}$, assuming
it holds for $\mathcal B_{k}$, $P_{k}$. Define 
\begin{equation}
  \label{eq:Xik}
\Xi_k(x,p):=\Phi_k(x,p)+\langle P_kx,x\rangle-2\langle r_k,x\rangle+\alpha_k
\end{equation}
Then $\mathcal 
B_{k+1}(p)=\min_{x_k} \Xi_k(x_k,p)$ due to \eqref{eq:Bdfn}. Combining $P_k\ge 0$ with definition of $\Phi_k$, we derive\footnote{$x\mapsto
  \langle Ax,x\rangle-2\langle x,q\rangle+c$ is convex iff $A$ is a symmetric  non-negative matrix.}  $x\mapsto\Xi_k(x,p)$ is a convex quadratic function for any $p$. This and $\Xi_k(x_k,p)\ge B_{k+1}(p)\ge0$ imply~\cite[p.268]{Rockafellar1970}
 $\mathrm{Argmin}_x\,\Xi_k(x,p)\ne\varnothing$. On the other hand~\cite[T.27.4]{Rockafellar1970}
$x\in\mathrm{Argmin}_{x}\,\Xi_k(x,p)$ iff $\nabla_{x}\Xi_k(x,p)=0$. Finally, we obtain $
\mathrm{Argmin}_x\,\Xi_k(x,p)\ne\varnothing$ and 
\begin{equation}
  \label{eq:XikArg}
\mathrm{Argmin}_{x}\,\Xi_k(x,p)=\{x:B_kx=C'_kS_kF_{k+1}p+r_k\}
\end{equation}
If we set $q_k=B^+_k(C'_kS_kF_{k+1}p+r_k)$ then $q_k\in\mathrm{Argmin}_x\,\Xi_k(x,p)$ due to~\cite{Albert1972}. 
Now, it is sufficient to calculate
$\Xi_k(q_k,p)$ in order to see that~\eqref{eq:BPk} holds for $\mathcal B_{k+1}$ and $P_{k+1}$. Assertion $P_{k+1}\ge0$ holds since
$\mathcal B_{k+1}$ is convex. To conclude the proof, let us define $\tilde x_\tau=p$ and $\tilde x_k\in\mathrm{Argmin}_x\,\Xi_k(x,\tilde x_{k+1})$, $k=\overline{1,\tau-1}$. Then $V_\tau(\tilde x_1\dots\tilde x_{\tau-1},p)=\mathcal B_\tau (p)$ due to \eqref{eq:Bdfn}-\eqref{eq:XikArg}. 
\end{pf}
\textbf{3.1. Example.} 
Consider a system $p_{k+1}=A_kp_k+v_k,p_0=v$ and assume 
$y_k=\tilde H_kp_k+g_k$ provided $p_k\in\mathbb R^n$ and $(v,g_0,\dots,g_\tau)$ belong to some ellipsoid. Now, given $y^*_1,\dots,y_\tau^*$, one needs to build the worst-case  estimation of $p^*_\tau$. We cannot apply directly standard minimax framework~\cite{Nakonechnii1978,Bertsekas1971,Baras1995} in this case as we do not have any information about the bounding set for $(v_0,\dots,v_\tau)$. Instead, we apply the approach\footnote{Since $\mathop{\mathrm{rank}}\bigl[
  \begin{smallmatrix}
    F_k\\H_k
  \end{smallmatrix}\bigr]< 2n$, it follows that results of~\cite{ishixara2} are not applicable.}, proposed above. Define $F_{k}:=(E,0)$, $C_k:=(A_k,E)$, $H_k:=(\tilde H_k,0)$ and $x^*_k:=(p^*_k,v^*_k)$. Then $x^*_k$ verifies ~\eqref{eq:state}-\eqref{eq:mes} with $f_0=v$, $f_k=0$ and $g_k=g^*_k$, $k=\overline{1,\tau}$. Therefore, the original problem may be reformulated as: given $y^*_1,\dots,y_\tau^*$, one needs to build the $l=(\ell,0)$-minimax estimation of $x^*_\tau$. Of course, the estimation of $x_\tau{}$ in the direction $l=(0,\ell)$ has an infinite minimax error for any $\ell$. But this is natural as $\langle l,x_\tau{}\rangle=\langle v_\tau{},\ell\rangle$ and $(v_1,\dots,v_\tau{})$ is unknown. \\
In what follows we present a numerical example. Let $p_k\in\mathbb R^2$, $\tilde H_k=
(\begin{smallmatrix}
 1&&0 
\end{smallmatrix})$, $A_k=\bigl(
\begin{smallmatrix}
  \frac 1{10}&&-\frac 1{5}\\\frac 7{25}&&-\frac 1{10}
\end{smallmatrix}
\bigr)$ and $v_k=(-\frac{k\sin(k)}{10},\frac{k\sin(k)}{10})$. We have generated $p^*_k$ with $v=(\frac 1{10},\frac 1{10})$ and $y^*_k$ with $g_k=\frac{2\sin(k)}{k+1}$, $R_k=\frac k{k+1}$, $S_k=\mathrm{diag}(1,1)$, $k=\overline{0,50}$. The results are displayed on Fig.~\ref{fig:2}.
\newline\textbf{3.2. Minimax estimator and $H_2/H_\infty$ filters.}
In \cite{Baras1995} a connection between set-membership state estimation
and $H_\infty$ approach is described for linear causal DAEs. The authors note
that the notion of informational state ($X(\tau)$ in our notation) is shown to
be intrinsic for both approaches: mathematical relations between
informational states of $H_\infty$ and set-membership state estimation are
described in \cite[Lemma 6.2.]{Baras1995}. Comparisons of set-membership
estimators with $H_2/H_\infty$ filters for linear
DAEs are presented in \cite{Sayed2001}, provided $F_k\equiv E$. Let us consider connections to $H_2$-filters in details. In~\cite{ishixara2} 
  the authors derive the Kalman's
  recursion to DAE from a deterministic least
  square fitting problem. Assuming $\mathop{\mathrm{rank}}\bigl[
  \begin{smallmatrix}
    F_k\\H_k
  \end{smallmatrix}\bigr]\equiv n$, they prove that the optimal estimation $\hat{x}_{i|k}$ 
can be found from 
  \begin{equation*}
    \begin{split}
     &\hat{x}_{k|k}=
   P_{k|k}F'_{k}A_{k-1}C_{k-1}\hat{x}_{k-1|k-1}+P_{k|k}H'_{k} R_{k}y_{k},\\
   &\hat{x}_{0|0}=P_{0|0}H'_0R_0y_0,
   A^{-1}_k=S^{-1}_k+C_{k}P_{k|k}C'_{k}\\
  &P^{-1}_{k|k}=F'_kA_{k-1}F_k+H'_kR_kH_k,
  P^{-1}_{0|0}=F'_0SF_0+H'_0R_0H_0
    \end{split}
  \end{equation*}
\begin{cor}
\label{l:eqfltr}  Let $r_0=H'_0R_0y_0$. If $\mathop{\mathrm{rank}}
  \bigl[\begin{smallmatrix}
    F_k\\H_k
  \end{smallmatrix}\bigr]\equiv n$  
then $I_k=0$ and $P^+_kr_k=\hat{x}_{k|k}$.
\end{cor}
\begin{pf}
Let us set $R_k=E,S=E,S_k=E$ for simplicity. 
The proof is by induction on $k$. For $k=0$, $P_{0|0}=P_0^{-1}$. The induction
hypothesis is $P_{k-1|k-1}=P_{k-1}^{-1}$. Suppose  
$A\in\mathbb R^{m\times n}$, $B\in\mathbb R^{n\times n}$, $B=B'>0$; then  
  \begin{equation}
    \label{eq:AS}
    A(A'A+B^{-1})^{-1}=(E+ABA')^{-1}AB
  \end{equation}
Using~\eqref{eq:AS} 
we get $
ABA'=[E+ABA']A[A'A+B^{-1}]^{-1}A'  
$. Combining this 
with the induction assumption we get $
E+C_{k-1}P_{k-1|k-1}C'_{k-1}=
    E+[E+C_{k-1}P_{k-1|k-1}C'_{k-1}]
C_{k-1}[P_{k-1}+C'_{k-1}C_{k-1}]^{-1}C'_{k-1}
$. 
By simple calculation it follows from the previous 
equality that 
    $E-C_{k-1}(P_{k-1}+C'_{k-1}C_{k-1})^{-1}C'_{k-1}=    
    (E+C_{k-1}P_{k-1|k-1}C'_{k-1})^{-1}$
Using this and definitions of $P_k$, $P_{k|k}$,
we get $P_k^{-1}=P_{k|k}$. 
\newline 
$P^{-1}_0r_0=\hat{x}_{0|0}$ due to corollary assumption. 
Suppose that $ 
P^{-1}_{k-1}r_{k-1}=\hat{x}_{k-1|k-1}$. 
The induction hypothesis, equality $P_k^{-1}=P_{k|k}$ and~\eqref{eq:AS} 
imply $
    (E+C_{k-1}P_{k-1|k-1}C'_{k-1})^{-1}C_{k-1}\hat{x}_{k-1|k-1}=
    C_{k-1}(C'_{k-1}C_{k-1}+P_{k-1})^{-1}_{k-1}r_{k-1}
$. 
Combining this with definitions of $\hat{x}_{k|k}$, $r_k$ we obtain  $
\hat{x}_{k|k}=P_k^{-1}r_k$. 
This concludes the proof.
\end{pf}
\section{Conclusion}\label{s:cnclsn}
We describe a set-membership state estimation approach for a linear operator
equation with uncertain disturbance restricted to belong to a convex
bounded closed subset of abstract Hilbert space. It is
based on the notion of an a posteriori set~\cite{Nakonechnii1978} $\mathscr
G(y)$, 
informational set \cite{Baras1995} and the notion of the minimax observable
subspace for the pair $(L,H)$. The latter is new for the set-membership state
estimation framework. It leads to nontrivial new 
results in set-membership state estimation for linear non-causal DAEs: we
present new equations describing the dynamics of the 
 minimax recursive estimator for discrete-time non-causal DAEs. 
We prove that these equations are consistent with the main results already
 established for regular DAEs. We illustrate benefits of considering 
non-causality in the state equation, applying our approach to a linear
filtration problem with unbounded noise. 
\begin{ack}
  It is a pleasure to thank Prof. A.Nakonechny and Dr. V.Pichkur for
  insightful discussions. 
  I also thank to Drs. V.Mallet, J.A. Hosking and my anonymous referees for their help with improving the
  presentation of this paper.   
\end{ack}
\bibliographystyle{plain}
\bibliography{lit}
\end{document}